# A THEOREM ON MAJORIZING MEASURES

By Witold Bednorz[1]

*Warsaw University*

Let $(T, d)$ be a metric space and $\varphi : \mathbb{R}_+ \to \mathbb{R}$ an increasing, convex function with $\varphi(0) = 0$. We prove that if $m$ is a probability measure $m$ on $T$ which is majorizing with respect to $d, \varphi$, that is, $\mathcal{S} := \sup_{x \in T} \int_0^{D(T)} \varphi^{-1}(\frac{1}{m(B(x,\varepsilon))}) \, d\varepsilon < \infty$, then

$$\mathbf{E} \sup_{s,t \in T} |X(s) - X(t)| \leq 32 \mathcal{S}$$

for each separable stochastic process $X(t)$, $t \in T$, which satisfies $\mathbf{E}\varphi(\frac{|X(s)-X(t)|}{d(s,t)}) \leq 1$ for all $s, t \in T$, $s \neq t$. This is a strengthening of one of the main results from Talagrand [*Ann. Probab.* **18** (1990) 1–49], and its proof is significantly simpler.

**1. Introduction.** In this paper, $(T, d)$ is a fixed metric space and $m$ a fixed probability measure (defined on Borel subsets) on $T$. We assume that $\text{supp}(m) = T$. For $x \in T$ and $\varepsilon \geq 0$, $B(x, \varepsilon)$ denotes the closed ball with center at $x$ and radius $\varepsilon$ [i.e., $B(x, \varepsilon) = \{y \in T : d(x, y) \leq \varepsilon\}$]. Let $D(T)$ be the diameter of $T$, that is, $D(T) = \sup\{d(s, t) : s, t \in T\}$. We define $C(T)$ as to be the space of all continuous functions on $T$ and $\mathcal{B}(T)$ as to be the space of all Borel and bounded functions on $T$.

For $a, b \geq 0$ we denote by $\mathcal{G}_{a,b}$ the class of all functions $\varphi : \mathbb{R}_+ \to \mathbb{R}$ which are increasing, continuous, which satisfy $\varphi(0) = 0$ and such that

$$(1.1) \qquad x \leq a + b \frac{\varphi(xy)}{\varphi(y)} \qquad \text{for all } x \geq 0, y \geq \varphi^{-1}(1).$$

For a fixed function $\varphi \in \mathcal{G}_{a,b}$ we define

$$\sigma(x) := \int_0^{D(T)} \varphi^{-1}\left(\frac{1}{m(B(x, \varepsilon))}\right) d\varepsilon,$$

Received February 2005; revised October 2005.
[1]Supported by Polish KBN-Grant 2 P03A 02722.

$$\bar{\mathcal{S}} := \int_T \sigma(u) m(du),$$

$$\mathcal{S} := \sup_{x \in T} \sigma(x).$$

We say that $m$ is a *majorizing measure* if $\mathcal{S} < \infty$. In the sequel we will use the convention that $0/0 = 0$.

The following theorem is the main result of the paper:

THEOREM 1.1. *If $\varphi$ is a Young function and $m$ is a majorizing measure on $T$, then, for each separable stochastic process $X(t)$, $t \in T$, such that*

(1.2) $$\sup_{s,t \in T} \mathbf{E}\varphi\left(\frac{|X(s) - X(t)|}{d(s,t)}\right) \leq 1,$$

*the following inequality holds:*

$$\mathbf{E} \sup_{s,t \in T} |X(s) - X(t)| \leq 32\mathcal{S}.$$

This is a generalization of Theorem 4.6 from Talagrand [3]. The method we use in this paper is new and the proof is simpler. Contrary to Talagrand's result, it works for all Young functions $\varphi$, in particular for $\varphi(x) \equiv x$. The author arrived at the idea of chaining with balls of given measure by studying [4] (see also [5]).

Our main tool needed to obtain Theorem 1.1 will be a Sobolev-type inequality.

THEOREM 1.2. *Suppose $\varphi \in \mathcal{G}_{a,b}$ and $R \geq 2$. Then there exists a probability measure $\nu$ on $T \times T$ such that, for each bounded, continuous function $f$ on $T$, the inequality*

$$\left|f(t) - \int_T f(u) m(du)\right| \leq aA\sigma(t) + bB\bar{\mathcal{S}} \int_{T \times T} \varphi\left(\frac{|f(u) - f(v)|}{d(u,v)}\right) \nu(du, dv),$$

*holds for all $t \in T$, where $A = \frac{R^3}{(R-1)(R-2)}$, $B = \frac{R^2}{R-1}$.*

An immediate consequence of Theorem 1.2 is the following corollary:

COROLLARY 1.1. *If $\varphi \in \mathcal{G}_{a,b}$ and $R \geq 2$ then there exists a probability measure $\nu$ on $T \times T$ such that, for all $f \in C(T)$,*

$$\sup_{s,t \in T} |f(s) - f(t)| \leq 2aA\mathcal{S} + 2bB\bar{\mathcal{S}} \int_{T \times T} \varphi\left(\frac{|f(u) - f(v)|}{d(u,v)}\right) \nu(du, dv),$$

*where $A = \frac{R^3}{(R-1)(R-2)}$, $B = \frac{R^2}{R-1}$.*

REMARK 1.1. In terms of absolutely summing operators, Corollary 1.1 means that the embedding of the Banach space of Lipschitz functions on $T$ into the Banach space of continuous and bounded functions on $T$ is $\varphi$-absolutely summing, as defined by Assouad [1].

Each increasing, convex function $\varphi$ with $\varphi(0) = 0$ (Young function) is in $\mathcal{G}_{1,1}$. Choosing $R = 4$, $a = b = 1$, Corollary 1.1 yields the following:

COROLLARY 1.2.  *If $\varphi$ is a Young function then there exists a probability measure $\nu$ on $T \times T$ such that, for all $f \in C(T)$,*

$$\sup_{s,t\in T} |f(s) - f(t)| \leq 32\mathcal{S}\left(\frac{2}{3} + \frac{1}{3}\int_{T\times T} \varphi\left(\frac{|f(u) - f(v)|}{d(u,v)}\right)\nu(du,dv)\right).$$

REMARK 1.2. For a Young function, it is usually possible to choose better constants than $a = b = 1$. For example, the function $\varphi(x) \equiv x$ is in $\mathcal{G}_{0,1}$. Setting $R = 2$, $a = 0$, $b = 1$ in Corollary 1.1, we obtain that there exists a probability measure $\nu$ on $T \times T$ such that

$$\sup_{s,t\in T} |f(s) - f(t)| \leq 8\bar{\mathcal{S}} \int_{T\times T} \frac{|f(u) - f(v)|}{d(u,v)}\nu(du,dv) \qquad \text{for all } f \in C(T).$$

The result is of interest if $\bar{\mathcal{S}} < \infty$, which is valid for a larger class of measures than majorizing measures.

We use Corollary 1.2 to prove the main result (Theorem 1.1).

## 2. Proofs and generalizations.

PROOF OF THEOREM 1.2. We can assume that $D(T) < \infty$, otherwise $\sigma(x) = \infty$, for all $x \in T$ and there is nothing to prove. There exists $k_0 \in \mathbb{Z}$ such that

$$R^{k_0} \leq \varphi^{-1}(1) < R^{k_0+1}.$$

For $x \in T$ and $k > k_0$ we define

(2.1) $$r_k(x) := \min\left\{\varepsilon \geq 0 : \varphi^{-1}\left(\frac{1}{m(B(x,\varepsilon))}\right) \leq R^k\right\}.$$

If $k = k_0$, we put $r_{k_0}(x) := D(T)$.

LEMMA 2.1.  *For $k \geq k_0$, functions $r_k$ are 1-Lipschitz.*

PROOF. Indeed, $r_{k_0}$ is constant, and if $k > k_0$ then for each $s, t \in T$ we obtain from the definition
$$\varphi^{-1}\left(\frac{1}{m(B(s, r_k(t) + d(s,t)))}\right) \leq \varphi^{-1}\left(\frac{1}{m(B(t, r_k(t)))}\right) \leq R^k.$$
Hence $r_k(s) \leq r_k(t) + d(s,t)$ and similarly $r_k(t) \leq r_k(s) + d(s,t)$, which means $r_k$ is 1-Lipschitz. □

We have
$$\sum_{k \geq k_0} r_k(x)(R^k - R^{k-1})$$
$$\leq \sum_{k \geq k_0} (r_k(x) - r_{k+1}(x))R^k + \limsup_{k \to \infty} r_{k+1}(x) R^{k+1}$$
$$\leq \sum_{k \geq k_0} \int_{r_{k+1}(x)}^{r_k(x)} \varphi^{-1}\left(\frac{1}{m(B(x,\varepsilon))}\right) d\varepsilon$$
$$+ \limsup_{k \to \infty} \int_0^{r_{k+1}(x)} \varphi^{-1}\left(\frac{1}{m(B(x,\varepsilon))}\right) d\varepsilon$$
$$= \int_0^{D(T)} \varphi^{-1}\left(\frac{1}{m(B(x,\varepsilon))}\right) d\varepsilon.$$

Consequently,

(2.2) $$\sum_{k \geq k_0} r_k(x) R^k \leq \frac{R}{R-1} \sigma(x).$$

Let us denote $B_k(x) := B(x, r_k(x))$.

For each $k \geq k_0$, we define the linear operator $S_k : \mathcal{B}(T) \to \mathcal{B}(T)$ by the formula
$$S_k f(x) := \fint_{B_k(x)} f(u) m(du) := \frac{1}{m(B_k(x))} \int_{B_k(x)} f(u) m(du).$$

If $f, g \in \mathcal{B}(T)$, $k \geq k_0$, we can easily check that:

1. $S_k 1 = 1$;
2. if $f \leq g$ then $S_k f \leq S_k g$, hence $|S_k f| \leq S_k |f|$;
3. $S_{k_0} f = \int_T f(u) m(du)$ and hence $S_k S_{k_0} f = S_{k_0} f$;
4. if $f \in C(T)$ then $\lim_{k \to \infty} S_k f(x) = f(x)$.

The last property holds true since $\lim_{k \to \infty} r_k(x) = 0$.

LEMMA 2.2. *If $m > k \geq k_0$ then*

(2.3) $$S_m S_{m-1} \cdots S_{k+1} r_k \leq \sum_{i=k}^m 2^{i-k} r_i.$$

PROOF. First we will show that for $i, j \geq k_0$,

(2.4) $$S_i r_j \leq r_i + r_j.$$

Indeed, due to Lemma 2.1, we obtain $r_j(v) \leq r_i(u) + r_j(u)$ for each $v \in B_i(u) = B(u, r_i(u))$. Since $S_i r_j(u) = \fint_{B_i(u)} r_j(v) m(dv)$, it implies (2.4).

We will prove Lemma 2.2 by induction on $m$. For $m = k+1$, inequality (2.3) has the form $S_{k+1} r_k \leq r_k + 2 r_{k+1}$, and it follows by (2.4). Suppose that, for $m-1$ such that $m - 1 > k \geq k_0$, it is

$$S_{m-1} S_{m-2} \cdots S_{k+1} r_k \leq \sum_{i=k}^{m-1} 2^{i-k} r_i.$$

Applying (2.4) to the above inequality, we get

$$S_m S_{m-1} \cdots S_{k+1} r_k \leq S_m \sum_{i=k}^{m-1} 2^{i-k} r_i \leq \sum_{i=k}^{m-1} 2^{i-k}(r_i + r_m) \leq \sum_{i=k}^{m} 2^{i-k} r_i. \quad \square$$

Observe that

(2.5)
$$\sum_{k=k_0}^{m-1} \left( \sum_{i=k}^{m} 2^{i-k} r_i \right) R^k = \sum_{k=k_0}^{m-1} \sum_{i=k}^{m} \left( \frac{2}{R} \right)^{i-k} r_i R^i$$
$$\leq \sum_{j=0}^{\infty} \left( \frac{2}{R} \right)^j \sum_{i=k_0}^{m} r_i R^i$$
$$\leq \frac{R}{R-2} \sum_{i=k_0}^{\infty} r_i R^i.$$

By the properties 1–4 of the operators $S_k$, $k \geq k_0$, we get

(2.6)
$$\left| f(t) - \int_T f(u) m(du) \right| = \lim_{m \to \infty} |S_m f - S_m S_{m-1} \cdots S_{k_0} f|(t)$$
$$= \lim_{m \to \infty} \left| \sum_{k=k_0}^{m-1} S_m \cdots S_{k+2} S_{k+1}(I - S_k) f \right|(t)$$
$$\leq \lim_{m \to \infty} \sum_{k=k_0}^{m-1} S_m \cdots S_{k+2} |S_{k+1}(I - S_k) f|(t).$$

We can easily check that

$$S_{k+1}(I - S_k) f(w) = \fint_{B_{k+1}(w)} \fint_{B_k(u)} (f(u) - f(v)) m(dv) m(du),$$

which gives

$$|S_{k+1}(I - S_k) f|(w) \leq \fint_{B_{k+1}(w)} \fint_{B_k(u)} |f(u) - f(v)| m(dv) m(du).$$

Condition (1.1) implies that, for $v \in B_k(u)$,

$$(2.7) \qquad \frac{|f(u) - f(v)|}{R^{k+1}d(u,v)} \leq a + \frac{b}{\varphi(R^{k+1})}\varphi\left(\frac{|f(u) - f(v)|}{d(u,v)}\right).$$

For each $v \in B_k(u)$, we have that $d(u,v) \leq r_k(u)$, and for $w \in T$ it is $m(B_{k+1}(w)) \geq \frac{1}{\varphi(R^{k+1})}$. Thus, for $v \in B_k(u)$, the following inequality holds:

$$|f(u) - f(v)| \leq ar_k(u)R^{k+1} + bm(B_{k+1}(w))r_k(u)R^{k+1}\varphi\left(\frac{|f(u) - f(v)|}{d(u,v)}\right).$$

Consequently,

$$|S_{k+1}(I - S_k)f|(w) \leq aR^{k+1}S_{k+1}r_k(w)$$
$$+ b\int_T r_k(u)R^{k+1}\fint_{B_k(u)} \varphi\left(\frac{|f(u) - f(v)|}{d(u,v)}\right)m(dv)m(du).$$

By Lemma 2.2, $S_m \cdots S_{k+2}S_{k+1}r_k \leq \sum_{i=k}^m 2^{i-k}r_i$, therefore,

$$S_m \cdots S_{k+2}|S_{k+1}(I - S_k)f|(t)$$
$$\leq aR\sum_{i=k}^m 2^{i-k}r_i(t)R^k$$
$$+ bR\int_T r_k(u)R^k \fint_{B_k(u)} \varphi\left(\frac{|f(u) - f(v)|}{d(u,v)}\right)m(dv)m(du).$$

Using (2.5), (2.6) and then (2.2) we obtain

$$\left|f(t) - \int_T f(u)m(du)\right|$$
$$\leq a\frac{R^2}{R-2}\sum_{k=k_0}^\infty r_k(t)R^k$$
$$+ bR\sum_{k=k_0}^\infty \int_T r_k(u)R^k \fint_{B_k(u)} \varphi\left(\frac{|f(u) - f(v)|}{d(u,v)}\right)m(dv)m(du)$$
$$\leq aA\sigma(t) + bR\sum_{k=k_0}^\infty \int_T r_k(u)R^k \fint_{B_k(u)} \varphi\left(\frac{|f(u) - f(v)|}{d(u,v)}\right)m(dv)m(du),$$

where $A = \frac{R^3}{(R-1)(R-2)}$. Let $\nu$ be a probability measure on $T \times T$ defined by

$$\nu(g) := \frac{1}{M}\sum_{k=k_0}^\infty \int_T r_k(u)R^k \fint_{B_k(u)} g(u,v)m(dv)m(du) \qquad \text{for } g \in \mathcal{B}(T \times T),$$

where $M = \sum_{k=k_0}^{\infty} \int_T r_k(u) R^k m(du)$. By (2.2) we obtain an inequality $M \leq \frac{R}{R-1} \int_T \sigma(u) m(u) = \frac{R}{R-1} \bar{\mathcal{S}}$ and thus

$$\left| f(t) - \int_T f(u) m(du) \right| \leq aA\sigma(t) + bB\bar{\mathcal{S}} \int_{T \times T} \varphi\left(\frac{|f(u) - f(v)|}{d(u,v)}\right) \nu(du, dv),$$

where $B = \frac{R^2}{R-1}$. Theorem 1.2 is proved. □

There is a standard way to strengthen the obtained inequalities. We provide it here for the sake of completeness:

THEOREM 2.1. *Let $\psi: \mathbb{R}_+ \to \mathbb{R}$ be an increasing, continuous function with $\psi(0) = 0$, and $\alpha, \beta \geq 0$ such that*

(2.8) $$\psi(x) \leq \alpha + \beta \frac{\varphi(xy)}{\varphi(y)} \qquad \text{for all } x \geq 0, y \geq 0,$$

*where $\varphi \in \mathcal{G}_{a,b}$. Then, for each bounded, continuous functions $f$ on $T$, the following inequality holds:*

$$\sup_{t \in T} \psi\left(\frac{|f(t) - \int_T f(u) m(du)|}{K}\right)$$
$$\leq \alpha + \beta \int_{T \times T} \varphi\left(\frac{|f(u) - f(v)|}{d(u,v)}\right) \nu(du, dv),$$

*where $K = (aA + bB)\mathcal{S}$, and $A, B, \nu$ are as in Theorem 1.2.*

PROOF. Given function $f$, let $c$ be chosen in such a way that

$$\psi(c) = \alpha + \beta \int_{T \times T} \varphi\left(\frac{|f(u) - f(v)|}{d(u,v)}\right) \nu(du, dv).$$

By (2.8) we get, for all $u, v \in T$,

$$(\psi(c) - \alpha) \varphi\left(\frac{|f(u) - f(v)|}{cd(u,v)}\right) \leq \beta \varphi\left(\frac{|f(u) - f(v)|}{d(u,v)}\right).$$

Hence

$$\int_{T \times T} \varphi\left(\frac{|f(u) - f(v)|}{cd(u,v)}\right) \nu(du, dv)$$
$$\leq \frac{\beta}{\psi(c) - \alpha} \int_{T \times T} \varphi\left(\frac{|f(u) - f(v)|}{d(u,v)}\right) \nu(du, dv) = 1.$$

Therefore, by Theorem 1.2, we obtain

$$\frac{1}{c} \sup_{t \in T} \left| f(t) - \int_T f(u) m(du) \right|$$

$$\leq aA\sigma(t) + bB\bar{S} \int_{T\times T} \varphi\Big(\frac{|f(u)-f(v)|}{cd(u,v)}\Big)\nu(du,dv)$$
$$\leq (aA+bB)\mathcal{S} = K,$$

which is the same as $\sup_{t\in T}\frac{|f(t)-\int_T f(u)m(du)|}{K} \leq c$. Since $\psi$ is increasing, we get

$$\sup_{t\in T}\psi\Big(\frac{|f(t)-\int_T f(u)m(du)|}{K}\Big)$$
$$= \psi\Big(\sup_{t\in T}\frac{|f(t)-\int_T f(u)m(du)|}{K}\Big) \leq \psi(c)$$
$$= \alpha + \beta \int_{T\times T}\varphi\Big(\frac{|f(u)-f(v)|}{d(u,v)}\Big)\nu(du,dv). \qquad \Box$$

REMARK 2.1. Similarly, we can prove that, for each $f \in C(T)$, the following inequality holds:

$$\sup_{s,t\in T}\psi\Big(\frac{|f(s)-f(t)|}{2K}\Big) \leq \alpha + \beta\int_{T\times T}\varphi\Big(\frac{|f(u)-f(v)|}{d(u,v)}\Big)\nu(du,dv).$$

Each Young function satisfies (1.1) with $a=1$, $b=1$. The minimal constant $K = (A+B)\mathcal{S} = \frac{2R^2}{R-2}\mathcal{S}$ is equal to $16\mathcal{S}$ and is attained for $R=4$. Let us consider functions $\varphi_p(x) \equiv x^p$, $p \geq 1$. The condition (1.1) is satisfied if and only if $(aq)^{1/q}(bp)^{1/p} \geq 1$, where $q = \frac{p}{p-1}$. Elementary calculations show that by choosing

$$R_p = 2 + \frac{1}{q}\Big(\Big(3q-\frac{q}{p}\Big)^{1/2}+1\Big),$$
$$a_p = \frac{1}{q}\Big(3q-\frac{q}{p}\Big)^{-1/(2p)},$$
$$b_p = \frac{1}{p}\Big(3q-\frac{q}{p}\Big)^{1/(2q)},$$

we obtain the minimal constant $K_p := 2(\frac{3p-1}{p})(3q-\frac{q}{p})^{1/(2q)}\mathcal{S}$.

Since $\varphi_p(x) \equiv x^p$ satisfies (2.8) for $\alpha=0, \beta=1$, we can conclude the above considerations with the following proposition:

PROPOSITION 2.1. *If $m$ is a majorizing measure on $T$, then there exists a probability measure $\nu$ on $T\times T$ such that*

$$\sup_{s,t\in T}|f(s)-f(t)|^p \leq (2K_p)^p\int_{T\times T}\Big(\frac{|f(u)-f(v)|}{d(u,v)}\Big)^p\nu(du,dv),$$

*for all $f \in C(T)$, where $K_p = 2(\frac{3p-1}{p})(3q-\frac{q}{p})^{1/(2q)}\mathcal{S}$.*

**3. An application to sample boundedness.** The theorems from the preceding section allow us to prove results concerning the boundedness of stochastic processes. In this paper we consider only separable processes. For such a process $X(t)$, $t \in T$, we have

$$\mathbf{E}\sup_{t \in T} X(t) := \sup_{F \subset T} \mathbf{E}\sup_{t \in F} X(t),$$

where the supremum is taken over all finite sets $F \subset T$.

THEOREM 3.1. *Suppose $\varphi \in \mathcal{G}_{a,b}$ is a Young function, and $R \geq 2$. For each process $X(t)$, $t \in T$, which satisfies (1.2), the following inequality holds:*

$$\mathbf{E}\sup_{s,t \in T} |X(s) - X(t)| \leq 2aA\mathcal{S} + 2bB\bar{\mathcal{S}},$$

*where $A = \frac{R^3}{(R-1)(R-2)}$, $B = \frac{R^2}{R-1}$.*

PROOF. Our argument follows the proof of Theorem 2.3, [3]. The process $X(t)$ $t \in T$, is defined on a probability space $(\Omega, \mathcal{F}, \mathbf{P})$. Take any point $t_0 \in T$. Condition (1.2) implies $\mathbf{E}|X(t) - X(t_0)| < \infty$, for all $t \in T$.

We define $Y(t) := X(t) - X(t_0)$. Necessarily, $\mathbf{E}|Y(t)| < \infty$, for all $t \in T$, condition (1.2) holds and $\mathbf{E}\sup_{s,t \in T}|X(s) - X(t)| = \mathbf{E}\sup_{s,t \in T}|Y(s) - Y(t)|$. First, we suppose that $\mathcal{F}$ is finite. We may identify points in each atom of $\mathcal{F}$, so we can assume that $\Omega$ is finite. Let us observe that

$$|Y(s,\omega) - Y(t,\omega)| \leq d(s,t)\varphi^{-1}(1/\mathbf{P}(\{\omega\})),$$

so trajectories of $Y$ are Lipschitz and consequently continuous. Using Corollary 1.1, the Fubini theorem and condition (1.2), we obtain

$$\mathbf{E}\sup_{s,t \in T} |Y(s) - Y(t)| \leq 2aA\mathcal{S} + 2bB\bar{\mathcal{S}} \int_{T \times T} \mathbf{E}\varphi\left(\frac{|Y(u) - Y(v)|}{d(u,v)}\right)\nu(du,dv)$$

$$= 2aA\mathcal{S} + 2bB\bar{\mathcal{S}}.$$

In the general case, we have to show that, for any finite $F \subset T$,

(3.1) $$\mathbf{E}\sup_{s,t \in F} |Y(s) - Y(t)| \leq 2aA\mathcal{S} + 2bB\bar{\mathcal{S}},$$

so we may assume that $\mathcal{F}$ is countably generated. There exists an increasing sequence $\mathcal{F}_n$ of finite $\sigma$-fields whose union generates $\mathcal{F}$. Since $\mathbf{E}|Y(t)| < \infty$, it is possible to define $Y_n(t) = \mathbf{E}(Y(t)|\mathcal{F}_n)$. Jensen's inequality shows that

$$\mathbf{E}\varphi\left(\frac{|Y_n(s) - Y_n(t)|}{d(s,t)}\right) \leq \mathbf{E}\varphi\left(\frac{|Y(s) - Y(t)|}{d(s,t)}\right) \leq 1.$$

We get (3.1) since $Y_n(t) \to Y(t)$, $\mathbf{P}$-a.s. and in $L_1$ for each $t \in F$. $\square$

Each Young function $\varphi \in \mathcal{G}_{1,1}$ and $\bar{\mathcal{S}} \leq \mathcal{S}$, so choosing $R = 4$, $a = b = 1$ in Theorem 3.1, we obtain Theorem 1.1.

REMARK 3.1. Our assumption that $\varphi$ is a Young function is not necessary. Suppose we have an arbitrary function $\varphi \in \mathcal{G}_{a,b}$ and $R \geq 2$. For each process $X(t)$, $t \in T$ which satisfies (1.2), the following inequality holds:

$$\mathbf{E} \sup_{s,t \in T} |X(s) - X(t)| \leq 4K,$$

where $K = (aA + bB)\mathcal{S}$, $A = \frac{R^3}{(R-1)(R-2)}$, $B = \frac{R^2}{R-1}$.

PROOF. Following the proof of Theorem 11.9 from [2], for every finite $F \subset T$, there exists a measurable map $f : T \to F$ such that $d(f(t), x) \leq 2d(t, x)$, for all $t \in T, x \in F$.

We define $\mu_F = f(m)$ so that $\mu_F$ is supported by $F$. Thus, $f(B(x, \varepsilon)) \subset B_F(x, 2\varepsilon)$, and finally we get $m(B(x, \varepsilon)) \leq \mu_F(B_F(x, 2\varepsilon))$. Since the process $X$ is continuous on $F$, similarly as in the proof of Theorem 3.1, we get

$$\mathbf{E} \sup_{s,t \in F} |X(s) - X(t)|$$
$$\leq 2(aA + bB) \sup_{x \in F} \int_0^{D(F)} \varphi^{-1}\left(\frac{1}{\mu_F(B(x, \varepsilon))}\right) d\varepsilon$$
$$\leq 2(aA + bB) \sup_{x \in F} \int_0^{D(F)} \varphi^{-1}\left(\frac{1}{m(B(x, 1/2\varepsilon))}\right) d\varepsilon \leq 4K. \quad \square$$

The method presented in Theorem 2.1 allows us to obtain the following result:

THEOREM 3.2. Let $\varphi, \psi$ be as in Theorem 2.1. For each process which satisfies (1.2), the following inequality holds:

$$\mathbf{E} \sup_{s,t \in T} \psi\left(\frac{|X(s) - X(t)|}{2K}\right) \leq \alpha + \beta,$$

where $K = (aA + bB)\mathcal{S}$, $A = \frac{R^3}{(R-1)(R-2)}$, $B = \frac{R^2}{R-1}$.

REMARK 3.2. In the case of function $\varphi_p(x) = x^p$, $p \geq 1$, following Remark 2.1, we obtain

$$\left\| \sup_{s,t \in T} |X(s) - X(t)| \right\|_p \leq 2K_p.$$

**Acknowledgment.** I would like to thank professor M. Talagrand for all his helpful comments.

DEPARTMENT OF MATHEMATICS
WARSAW UNIVERSITY
UL. BANACHA 2
WARSAW 02-097
POLAND
E-MAIL: wbednorz@mimuw.edu.pl